\documentclass[preprint,12pt]{elsarticle}

\usepackage{amssymb,amsmath,verbatim,tabularx,enumitem,colonequals,graphicx,psfrag,subfigure,xcolor,tikz,rotating} 
\usepackage[amsmath,amsthm,thmmarks]{ntheorem}

\newtheorem{thm}{Theorem}
\newtheorem{prop}[thm]{Proposition}
\newtheorem{lem}[thm]{Lemma}
\newtheorem{cor}[thm]{Corollary}
\newtheorem{prob}[thm]{Problem}
\newdefinition{defn}{Definition}
\newtheorem{example}{Example}
\newproof{pf}{Proof}

\newcommand{\LB}{\textnormal{LB}}
\newcommand{\UB}{\textnormal{UB}}

\newcommand{\WS}{\textnormal{WS}}
\newcommand{\Budget}{\textnormal{Budget}}
\newcommand{\opt}{\textsc{opt}}
\newcommand{\alg}{\textsc{alg}}
\newcommand{\bgamma}{\bar{\gamma}}
\newcommand{\CLIQUE}{\textnormal{CLIQUE}}
\newcommand{\twoCLIQUE}{\textnormal{2-CLIQUE}}

\usepackage{array}
\newcolumntype{L}[1]{>{\raggedright\let\newline\\\arraybackslash\hspace{0pt}}m{#1}}
\newcolumntype{C}[1]{>{\centering\let\newline\\\arraybackslash\hspace{0pt}}m{#1}}
\newcolumntype{R}[1]{>{\raggedleft\let\newline\\\arraybackslash\hspace{0pt}}m{#1}}

\makeatletter
\def\ps@pprintTitle{%
   \let\@oddhead\@empty
   \let\@evenhead\@empty
   \def\@oddfoot{\reset@font\hfil\thepage\hfil}
   \let\@evenfoot\@oddfoot
}
\makeatother

\begin{document}
\begin{frontmatter}

\title{A General Approximation Method\\ for Bicriteria Minimization Problems}

\author[kob]{Pascal Halffmann}
\ead{halffmann@uni-koblenz.de}

\author[kob]{Stefan Ruzika}
\ead{ruzika@uni-koblenz.de}

\author[kl]{Clemens Thielen\corref{cor1}}
\ead{thielen@mathematik.uni-kl.de}

\author[kob]{David Willems}
\ead{davidwillems@uni-koblenz.de}

\cortext[cor1]{Corresponding author. Fax: +49 (631) 205-4737. Phone: +49 (631) 205-4590}

\address[kob]{Mathematical Institute, University of Koblenz-Landau,
Campus Koblenz, D-56070~Koblenz, Germany}

\address[kl]{Department of Mathematics, University of Kaiserslautern,
  Paul-Ehrlich-Str.~14, D-67663~Kaiserslautern, Germany}
    
\begin{abstract}
We present a general technique for approximating bicriteria minimization problems with positive-valued, polynomially computable objective functions. Given $0<\epsilon\leq1$ and a polynomial-time $\alpha$-approximation algorithm for the corresponding weighted sum problem, we show how to obtain a bicriteria $(\alpha\cdot(1+2\epsilon),\alpha\cdot(1+\frac{2}{\epsilon}))$-approximation algorithm for the budget-constrained problem whose running time is polynomial in the encoding length of the input and linear in $\frac{1}{\epsilon}$.

  Moreover, we show that our method can be extended to compute an $(\alpha\cdot(1+2\epsilon),\alpha\cdot(1+\frac{2}{\epsilon}))$-approximate Pareto curve under the same assumptions. Our technique applies to many minimization problems to which most previous algorithms for computing approximate Pareto curves cannot be applied because the corresponding gap problem is $\textsf{NP}$-hard to solve. For maximization problems, however, we show that approximation results similar to the ones presented here for minimization problems are impossible to obtain in polynomial time unless $\textsf{P}=\textsf{NP}$.
\end{abstract}

\begin{keyword}
	multicriteria optimization \sep bicriteria approximation algorithm \sep approximate Pareto curve
\end{keyword}
\end{frontmatter}

\section{Introduction}
Multicriteria optimization is one of the fastest growing fields of research in optimization and operations research. It provides methods and techniques for solving optimization problems with multiple, equally important but opposing objectives. In the last decades, a variety of theoretical results and algorithms have been developed to improve solvability of multicriteria optimization problems.

\medskip

An important class of problems that has recently received much attention from researchers are \emph{combinatorial} multicriteria optimization problems. Many of these problems, however, are already $\textsf{NP}$-hard to solve when only a single objective function is considered, so computing exact (Pareto or nondominated) solutions or even the complete Pareto curve for the multicriteria versions is often intractable. Even for problems where the unicriterion version is efficiently solvable, computing the complete Pareto curve may turn out to be intractable due to the possibly exponential number of nondominated solutions. This motivates to study approximations of the Pareto curve that can be computed in polynomial time.

\medskip

For a bicriteria optimization problem, a common alternative approach to considering Pareto solutions is to optimize one of the two objective functions subject to a bound (often called \emph{budget constraint}) on the value of the other objective function. Unfortunately, the budget-constrained problem often turns out to be $\textsf{NP}$-hard to solve or even approximate even if the corresponding unicriterion problem without the budget constraint is solvable in polynomial time. This motivates to study polynomial-time \emph{bicriteria approximation algorithms} for the budget-constrained problem. Such algorithms compute solutions that violate the budget constraint by at most a given factor and at the same time obtain a bounded approximation factor with respect to the objective function that has to be optimized.

\medskip

There are specialized bicriteria approximation algorithms for the budget-constrained versions of many particular bicriteria optimization problems and also many dedicated algorithms that compute approximations of the Pareto curve of particular bicriteria problems. However, besides the seminal work of Papadimitriou and Yannakakis~\cite{papadimitriou2000approximability} (and the extensions due to Vassilvitskii and Yannakakis~\cite{vassilvitskii2005efficiently} and Diakonikolas and Yannakakis~\cite{diakonikolas2009small}) and the results of Gla\ss er et al.~\cite{glasser2010approximability,Glasser+etal:multi-hardness}, there is no general approximation algorithm for the Pareto curve that applies to a broad class of bicriteria problems. Moreover, the algorithms based on the method of Papadimitriou and Yannakakis compute approximate Pareto curves in polynomial time only for those problems whose corresponding \emph{gap problem} can be solved in polynomial time, which is not possible for many important problems.

\medskip

In this article, we present a general bicriteria approximation algorithm that applies to the budget-constrained version of a broad class of bicriteria minimization problems. Moreover, we show that our technique also yields a polynomial-time algorithm for computing an approximate Pareto curve even for problems for which the gap problem is $\textsf{NP}$-hard to solve.

\subsection{Previous Work}
Papadimitriou and Yannakakis~\cite{papadimitriou2000approximability} present a method for generating an approximate Pareto curve obtaining an approximation factor of $1+\epsilon$ in every objective function (called an \emph{$\epsilon$-Pareto curve}) for general multicriteria minimization and maximization problems with positive-valued, polynomially computable objective functions. They show that an $\epsilon$-Pareto curve with size polynomial in the encoding length of the input and $\frac{1}{\epsilon}$ always exist; but the construction of such a curve is only possible in (fully) polynomial time if the following \emph{gap problem} (which we state for minimization problems here) can also be solved in (fully) polynomial time:

\begin{prob}[Gap Problem]
	Given an instance~$I$ of a $k$-criteria minimization problem, a vector~$b\in \mathbb{R}^k$, and $\epsilon>0$, either return a feasible solution~$s$ whose objective value~$f(s)$ satisfies $f_i(s)\leq b_i$ for all~$i$ or answer that there is no feasible solution~$s'$ with $f_i(s')\leq \frac{b_i}{1+\epsilon}$ for all~$i$.
\end{prob}

Here, as usual in the context of approximation, \emph{polynomial time} refers to a running time that is polynomial in the encoding length of the input and \emph{fully polynomial time} refers to a running time that is additionally polynomial in $\frac{1}{\epsilon}$. Hence, if the gap problem is solvable in (fully) polynomial time, the result of Papadimitriou and Yannakakis~\cite{papadimitriou2000approximability} shows that an $\epsilon$-Pareto curve can be constructed in (fully) polynomial time (i.\,e., the problem admits a multicriteria \emph{(fully) polynomial-time approximation scheme} abbreviated by \emph{MPTAS} (\emph{MFPTAS})).

\medskip

In a succeeding paper, Vassilvitskii and Yannakakis~\cite{vassilvitskii2005efficiently} show that, for bicriteria problems whose gap problem can be solved in (fully) polynomial time, an $\epsilon$-Pareto curve that has cardinality at most three times the cardinality of the smallest $\epsilon$-Pareto curve~$P^*_{\epsilon}$ can be constructed in (fully) polynomial time. Moreover, they showed that the factor of three is best possible in the sense that, for some problems, it is $\textsf{NP}$-hard to do better. In a more recent paper, Bazgan et al.~\cite{Bazgan+etal:min-pareto} present an algorithm that also works for bicriteria problems whose gap problem can be solved in (fully) polynomial time. This algorithm also computes an $\epsilon$-Pareto curve that has cardinality at most three times~$|P^*_{\epsilon}|$ in (fully) polynomial time, but uses a subroutine~$\textsf{SoftRestrict}$ for a different subproblem that is polynomially equivalent to the gap problem.

Diakonikolas and Yannakakis~\cite{diakonikolas2009small} show that, for bicriteria problems for which the budget-constrained problem (for some choice of the budgeted objective function) admits a (fully) polynomial-time approximation scheme, an $\epsilon$-Pareto curve that has at most twice as many elements as the smallest $\epsilon$-Pareto curve~$P^*_{\epsilon}$ can be constructed in (fully) polynomial time.
Again, this factor of two is shown to be best possible.

\medskip 

Besides the seminal work of Papadimitriou and Yannakakis and the extensions mentioned above, there is another general approximation method for bicriteria optimization problems under the assumption of positive-valued, polynomially computable objective functions. This can be obtained when combining two results of Gla\ss er et al.~\cite{glasser2010approximability,Glasser+etal:multi-hardness}. They introduce the following problem similar to the gap problem above:

\begin{prob}[Approximate Domination Problem]
	Given an instance~$I$ of a $k$-criteria minimization problem, a vector~$b\in \mathbb{R}^k$, and $\alpha\geq 1$, either return a feasible solution~$s$ whose objective value~$f(s)$ satisfies $f_i(s)\leq \alpha\cdot b_i$ for all~$i$ or answer that there is no feasible solution~$s'$ with $f_i(s')\leq b_i$ for all~$i$.
\end{prob}

Gla\ss er et al. show that, if this problem is solvable in polynomial time for some $\alpha\geq 1$, then an approximate Pareto curve obtaining an approximation factor of $\alpha\cdot (1+\epsilon)$ in every objective function can be computed in fully polynomial time.
Moreover, they show that the approximate domination problem with $\alpha\colonequals k\cdot\delta$ can be solved by using a $\delta$-approximation algorithm for the weighted sum problem of the $k$-criteria problem (i.e., for the unicriterion problem obtained by minimizing the sum of the $k$~objective functions, where each objective function is assigned a positive coefficient / weight).

Together, this implies that an approximate Pareto curve obtaining an approximation factor of $k\cdot\delta\cdot(1+\epsilon)$ in every objective function can be computed in fully polynomial time for $k$-criteria minimization problems provided that the objective functions are positive-valued and polynomially computable and a $\delta$-approximation algorithm for the weighted sum problem exists. As this result is not explicitly stated in~\cite{glasser2010approximability,Glasser+etal:multi-hardness}, no bounds on the running time are provided. 

\subsection{Our Contribution}
We consider general bicriteria minimization problems under the usual assumption that the objective functions are polynomially computable and take on only strictly positive values. Given $0<\epsilon\leq1$ and a polynomial-time $\alpha$-approxima\-tion algorithm for the corresponding weighted sum problem, we propose a general bicriteria $(\alpha\cdot(1+2\epsilon),\alpha\cdot(1+\frac{2}{\epsilon}))$-approximation algorithm for the budget-constrained problem whose running time is polynomial in the encoding length of the input and linear in $\frac{1}{\epsilon}$.
In the case that an \emph{exact} polynomial-time algorithm is given for the weighted sum problem, we show how the running time of our algorithm can be further improved by using binary search. If this exact algorithm for the weighted sum problem additionally satisfies the necessary assumptions for applying Megiddo's parametric search technique~\cite{Megiddo:rational-objectives,Megiddo:parallel}, we show that the approximation guarantee can be improved to $(1+\epsilon,1+\frac{1}{\epsilon})$ and the running time becomes strongly polynomial (provided that the weighted sum algorithm runs in strongly polynomial time).

\medskip

Moreover, we show that our algorithm is fit to compute an $(\alpha\cdot(1+2\epsilon),\alpha\cdot(1+\frac{2}{\epsilon}))$-approximate Pareto curve in time polynomial in the encoding length of the input and linear in $\frac{1}{\epsilon}$ under the same assumptions that are needed for the basic version of the bicriteria approximation algorithm.\footnote{In Appendix~$1$, we show that the assumption of strictly positive objective values can be relaxed slightly by only assuming nonnegativity of the objective values in the case that we can compute positive upper and lower bounds on all strictly positive objective values of feasible solutions in polynomial time.} Furthermore, if a \emph{parametric} $\alpha$-approximation algorithm for the weighted sum problem exists (which computes a sequence of $\alpha$-approximate solutions for the weighted sum problem for all combinations of positive weights simultaneously in polynomial time), we show that our method yields an $(\alpha\cdot(1+\epsilon),\alpha\cdot(1+\frac{1}{\epsilon}))$-approximate Pareto curve and the running time becomes strongly polynomial (again provided that the parametric approximation algorithm for the weighted sum problem runs in strongly polynomial time).

\medskip

Our method applies to a large number of bicriteria minimization problems. In particular, it can be employed on problems whose unicriterion version is $\textsf{APX}$-complete and coincides with the weighted sum problem. For these problems, the gap problem is not solvable in polynomial time unless $\textsf{P}=\textsf{NP}$, so the results of Papadimitriou and Yannakakis~\cite{papadimitriou2000approximability} and the succeeding papers~\cite{vassilvitskii2005efficiently,diakonikolas2009small,Bazgan+etal:min-pareto} cannot be used. We demonstrate examples of such problems and a comparison of our results for the specific problems to previous results in Section~\ref{sec:applications}.  

\medskip

Of course, a natural question is whether analogous approximation results can also be obtained for bicriteria \emph{maximization} problems. We show, however, that similar approximation results are impossible to obtain in polynomial time for general maximization problems unless $\textsf{P}=\textsf{NP}$.

\medskip

We remark that, in all versions of our method, every choice of $0<\epsilon\leq1$ yields a different trade-off in the approximation guarantees obtained with respect to the two objective functions. While the approximation guarantee obtained with respect to the first objective function decreases as $\epsilon$ is decreased, the approximation guarantee with respect to the second objective function \emph{increases} for smaller values of~$\epsilon$.
However, by exchanging the role of the two objective functions, all our results also hold with the approximation guarantees reversed, so the basic version of our method can also be used for obtaining an approximation guarantee of $(\alpha\cdot(1+\frac{2}{\epsilon}),\alpha\cdot(1+2\epsilon))$. 
Still, the behavior of the approximation guarantees in our method is different from many other bicriteria approximation techniques, where a decrease in the error parameter usually yields an improved approximation guarantee with respect to both objectives.

\medskip

The technique we use for proving the approximation guarantee of our algorithm is similar to the technique utilized by Eubank et al.~\cite{Eubank+etal:epidemics}, who present a
bicriteria $(1+2\epsilon,1+\frac{2}{\epsilon})$-approximation algorithm for the \emph{quarantining problem}. In this problem, one is given an undirected graph~$G=(V,E)$, an
initial infected set~$I\subseteq V$, a budget~$B$, and a cost~$c(e)$ for each edge~$e\in E$ and the task consists of computing a cut~$(S,V\setminus S)$ in~$G$ such that
$I\subseteq S$, $|S|\leq B$, and the cost of the edges crossing the cut is minimized. In their bicriteria approximation algorithm, Eubank et al.~\cite{Eubank+etal:epidemics}
extend the graph~$G$ by adding a source~$s$ and a sink~$t$ in a suitable way with the cost of some new edges given by a positive parameter~$\beta$. Then they show that, for a minimum $s$-$t$-cut~$(S,\bar{S})$ in the extended graph, the sum of $\beta\cdot |S\setminus I|$ and the cost of the original edges in the cut is at most~$B+\beta\cdot\opt$. When choosing $\beta\colonequals \frac{\epsilon\cdot B}{\opt}$, bounding each of the two summands by $B+\beta\cdot\opt$ yields an approximation guarantee of $(1+\epsilon,1+\frac{1}{\epsilon})$ and approximating the correct value of~$\beta$ by searching in multiplicative steps of~$(1+\epsilon)$ still yields an approximation guarantee of $(1+2\epsilon,1+\frac{2}{\epsilon})$. It should be noted that this argumentation is a special case of our arguments used to prove Propositions~\ref{prop:approx} and~\ref{prop:approx2}. However, our results are much more general since, by using the weighted sum problem instead of the graph extension argument, it applies to arbitrary bicriteria minimization problems with positive objective function values (in which no graph structure is needed and the set of solutions need not even be discrete). Moreover, we generalize the argumentation to the case that the auxiliary problem (in our case the weighted sum problem) can be solved only approximately, whereas an exact solution of the auxiliary minimum cut problem is needed in the argumentation in~\cite{Eubank+etal:epidemics}. In the case that an exact algorithm is available for the auxiliary problem, we improve the running time of the procedure drastically by using binary search and parametric search (and even obtain a strongly polynomial running time in some cases). Moreover, we show how an extension of our method can be used to compute approximate Pareto curves by applying the algorithm for the auxiliary problem to a suitably chosen range of parameters.

\medskip

We remark that, in a previous article~\cite{marathe1998networkdesign} by partially the same authors as~\cite{Eubank+etal:epidemics}, a technique different from the one used in~\cite{Eubank+etal:epidemics} has already been utilized in order to design a bicriteria $(\alpha\cdot(1+\epsilon),\alpha\cdot(1+\frac{1}{\epsilon}))$-approximation algorithm for the budget-constrained version of bicriteria network design problems from an $\alpha$-approximation algorithm for the corresponding weighted sum problem. However, one can show that the algorithm presented in~\cite{marathe1998networkdesign} is in general \emph{not} applicable unless an \emph{exact} algorithm is used for the weighted sum problem. Moreover, even if an exact algorithm is used for the weighted sum problem, the algorithm might fail to output a solution due to an ill-defined search interval even though the considered problem instance is feasible (we provide explicit examples for both phenomena in Appendix~$2$).

\medskip

The remainder of the paper is organized as follows: In Section~\ref{sec:prob-def}, we formally introduce the class of bicriteria minimization problems we consider and provide the necessary definitions concerning bicriteria approximation algorithms and (approximate) Pareto curves. Section~\ref{sec:bicriteria-approx} presents our general bicriteria approximation algorithm for the budget-constrained problem. In Section~\ref{sec:running-time-imp}, we show how the running time of the algorithm can be improved by using binary search or parametric search for problems for which the weighted sum problem can be solved exactly in polynomial time. Section~\ref{sec:Pareto-curve} presents the modification for computing an approximate Pareto curve and Section~\ref{sec:applications} provides applications of our technique to specific bicriteria minimization problems as well as a comparison of our results to previous results for these problems. In Section~\ref{sec:maximization}, we show that approximation results similar to the ones we obtain for minimization problems are not possible for general maximization problems unless $\textsf{P}=\textsf{NP}$. Section~\ref{sec:conclusion} summarizes our results and lists directions for future work.

\section{Problem Definition}\label{sec:prob-def}

We start by defining general bicriteria optimization problems and the different solution concepts considered within the paper. We denote by $\Pi$ a general bicriteria optimization problem. For a given instance~$I$ of~$\Pi$, we let $S(I)$ denote the \emph{set of feasible solutions} of instance~$I$. The two \emph{objective functions} of the problem that should be minimized are given by polynomial-time algorithms~$f_1,f_2$ that, given an instance~$I$ of $\Pi$ and a feasible solution~$x\in S(I)$, compute the values~$f_1(I,x),f_2(I,x)$, which are assumed to be positive rational numbers. We assume that, for $i=1,2$ and each fixed instance~$I$, there exist positive rational numbers~$\LB(I,i),\UB(I,i)$ with encoding length polynomial in the encoding length~$|I|$ of $I$ such that $\LB(I,i)\leq f_i(I,x) \leq \UB(I,i)$ for every feasible solution~$x\in S(I)$ (this is a consequence of the polynomial-time computability of the objective functions).

We let $Y(I)\colonequals \{f(I,x)=(f_1(I,x),f_2(I,x)): x\in S(I)\}$ denote the \emph{set of vectors of objective values} (or \emph{images}) of feasible solutions of instance~$I$. Usually, the considered instance~$I$ will be clear from context and we use the abbreviations $S\colonequals S(I)$, $Y\colonequals Y(I)$, $f_i(x)\colonequals f_i(I,x)$, $\LB(i)\colonequals\LB(I,i)$, and $\UB(i)\colonequals \UB(I,i)$ for $i=1,2$. We also refer to the objective value~$f_i(x)$ of a feasible solution~$x\in S$ with respect to the objective function~$f_i$ as the \emph{$f_i$-cost} of~$x$.

\begin{defn}\label{def:Pareto-curve}
For an instance~$I$ of $\Pi$, we say that an image~$y=f(I,x)\in Y(I)$ is \emph{dominated} by another image~$y'=f(I,x')\in Y(I)$ if $y'_i = f_i(I,x')\leq f_i(I,x) = y_i$
for $i=1,2$, but $y'\neq y$. If the image~$y=f(I,x)$ is not dominated by any other image~$y'$, we call~$y$ \emph{nondominated} and the feasible solution~$x\in S(I)$ \emph{efficient}. The set~$P(I)$ of all efficient solutions for instance~$I$ is called the \emph{Pareto curve} for instance~$I$.
\end{defn}

Note that the notion of the Pareto curve is somewhat blurred in the literature as it is sometimes used to describe the set of efficient solutions and sometimes to describe the set of nondominated images. We exclusively use the term Pareto curve to refer to the set of efficient solutions of a problem instance here in order to be consistent with the usual definition of approximate Pareto curves (see Definition~\ref{def:approx-Pareto}).

\medskip

One possibility to compute efficient solutions for a bicriteria optimization problem~$\Pi$ is to solve the \emph{weighted sum problem} associated with~$\Pi$ (introduced by Zadeh~\cite{zadeh1963optimality}):

\begin{defn}\label{def:weighted-sum}
For a given instance~$I$ of $\Pi$ and a positive weight~$\gamma>0$, the \emph{weighted sum problem}~$\Pi^{\WS}(\gamma)=\Pi^{\WS}(I,\gamma)$ asks for a feasible solution~$x\in S(I)$ that minimizes the weighted
sum $f_1(I,x) + \gamma\cdot f_2(I,x)$ over all feasible solutions in~$S(I)$.
\end{defn}

Note that, in the above definition, we have normalized the weight assigned to the first objective function~$f_1$ to~$1$, which is without loss of generality as long as strictly positive weights for both objective functions are considered. This is motivated by the well-known fact that, for strictly positive weights for both objective functions, every optimal solution of the weighted sum problem is an efficient solution for the corresponding instance of the bicriteria problem (efficient solutions that can be obtained in this way are called \emph{supported} efficient solutions; all other efficient solutions are called \emph{nonsupported} efficient solutions, cf.~\cite{Ehrgott:book}).

\medskip

Another possible way to look at a bicriteria optimization problem is to turn it into a unicriterion problem by imposing a bound on one of the two objectives and optimize the other:

\begin{defn}
For a bicriteria optimization problem~$\Pi$, the \emph{budget-con\-strained problem}~$\Pi^{\Budget}$ is the problem that, given an instance~$I$ of $\Pi$ and an upper bound~$B>0$ on the value of the first objective function~$f_1$, asks for a feasible solution~$x\in S(I)$ that minimizes the second objective function~$f_2$ over the set of all feasible solutions~$x\in S(I)$ for which
$f_1(I,x)\leq B$.
\end{defn}

Note that the budget-constrained problem is also often referred to as the \emph{$\epsilon$-constraint scalarization}, which was first considered by Haimes et al.~\cite{haimes1971bicriterion}. However, since the letter $\epsilon$ is used with a different meaning here, we exclusively refer to the problem as the budget-constrained problem throughout the rest of the paper.

\medskip

If a unicriterion problem such as the weighted sum problem~$\Pi^{\WS}(\gamma)$ or the budget-constrained problem~$\Pi^{\Budget}$ turns out to be $\textsf{NP}$-hard to solve, one is often interested in algorithms that compute approximate solutions in polynomial time:

\begin{defn}\label{def:approx-alg}
A \emph{(polynomial-time) $\alpha$-approximation algorithm} (with $\alpha\geq 1$) for a minimization problem is an algorithm that, for any given
instance~$I$ of encoding length~$|I|$, finds a feasible solution with objective value at most $\alpha$ times the optimal objective value in time bounded by a polynomial in~$|I|$ if the instance~$I$
admits a feasible solution, and outputs infeasibility of the instance after a number of steps bounded by a polynomial in $|I|$ if no feasible solution for instance~$I$ exists.
\end{defn}

A more general approximation approach for approximating problems such as the budget-constrained problem~$\Pi^{\Budget}$ is to study \emph{bicriteria approximation algorithms} that compute solutions that violate the budget constraint on the first objective function by at most a given factor and at the same time obtain a bounded approximation factor with respect to the second objective function:\footnote{
Algorithms violating a specific (budget) constraint of an optimization problem by a given factor are also often referred to as using \emph{resource augmentation}.}

\begin{defn}\label{def:bicriteria-approx}
A \emph{(polynomial-time) $(\alpha,\beta)$-approximation algorithm} (with $\alpha,\\ \beta \geq1$) for $\Pi^{\Budget}$ is an algorithm that, for any given instance~$I$ of $\Pi$ and any given budget~$B>0$, 
requires time polynomial in the encoding length of~$I, B$ to find a solution~$x\in S(I)$ such that the following holds: If there exists at least one feasible solution~$x'\in S(I)$ with $f_1(I,x')\leq B$,
then the solution~$x$ returned by the algorithm satisfies $f_1(I,x)\leq \alpha\cdot B$ and $f_2(I,x)\leq \beta\cdot \opt(B)$, where $\opt(B)\colonequals \inf\{f_2(I,\tilde{x}):\tilde{x}\in S(I), f_1(I,\tilde{x})\leq B\}\geq \LB(I,2)$ denotes the minimum $f_2$-cost among the feasible solutions with $f_1$-cost at most~$B$.
\end{defn}

Using a similar concept of approximation, one can study approximating the whole Pareto curve. Informally, an \emph{approximate Pareto curve} for an instance~$I$ of $\Pi$ is a set of feasible solutions for
instance~$I$ that approximately dominate all other solutions, which means that, for every feasible solution~$x\in S(I)$, the set contains a solution that has value at most a given factor larger
in each objective function.

\begin{defn}\label{def:approx-Pareto}
For an instance~$I$ of $\Pi$ and $\alpha,\beta\geq 1$, an $(\alpha,\beta)$-approximate Pareto curve is a set~$P_{\alpha,\beta}(I)$ of feasible solutions such that, for every feasible
solution~$x\in S(I)$, there exists a solution~$x'\in P_{\alpha,\beta}(I)$ with $f_1(I,x')\leq \alpha\cdot f_1(I,x)$ and $f_2(I,x')\leq \beta\cdot f_2(I,x)$. Moreover, for $\epsilon>0$, a $(1+\epsilon,1+\epsilon)$-approximate Pareto curve is called an \emph{$\epsilon$-Pareto curve}.
\end{defn}

\section{A General Bicriteria Approximation Algorithm}\label{sec:bicriteria-approx}
Given $0<\epsilon\leq 1$ and a polynomial-time $\alpha$-approximation algorithm~$\alg$ for the weighted sum problem~$\Pi^{\WS}$, we construct a bicriteria $(\alpha\cdot(1+2\epsilon),\alpha\cdot(1+\frac{2}{\epsilon}))$-approximation algorithm for $\Pi^{\Budget}$ whose running time is polynomial in the encoding length of~$I,B$ and in $\frac{1}{\epsilon}$.

\medskip

In the rest of this section, we consider a fixed instance~$I$ of~$\Pi$ and omit the explicit reference to the instance. We assume that we are given a budget~$B>0$ such that the instance of $\Pi^{\Budget}$ defined by~$I$ and~$B$ admits at least one feasible solution and let
\begin{align*}
\opt(B)\colonequals \inf\{f_2(x):x\in S, f_1(x)\leq B\}\geq \LB(2)>0
\end{align*}
denote the minimum $f_2$-cost among the feasible solutions for this instance.
Moreover, we let $x^{\alg}=x^{\alg}(\gamma)$ denote the solution returned by $\alg$ when applied to the instance of $\Pi^{\WS}$ defined by a weight~$\gamma>0$. The following two propositions provide the foundation of our algorithm:

\begin{prop}\label{prop:approx}
If $\gamma=\bgamma\colonequals\frac{\epsilon\cdot B}{\opt(B)}$, then $x^{\alg}=x^{\alg}(\gamma)$ satisfies $f_1(x^{\alg})\leq \alpha\cdot(1+\epsilon)\cdot B$ and
$f_2(x^{\alg})\leq\alpha\cdot\left(1+\frac{1}{\epsilon}\right)\cdot\opt(B)$.
\end{prop}

\begin{proof}
By definition of $\opt(B)$, we know that, for every $\mu>0$, there exists a feasible solution~$x^{\mu}\in S$ with $f_1(x^{\mu})\leq B$ and $f_2(x^{\mu})\leq \opt(B)+\mu$. The objective value of such a
solution~$x^{\mu}$ in the weighted sum problem~$\Pi^{\WS}$ with weight~$\gamma$ then satisfies
\begin{align*}
  f_1(x^{\mu}) + \gamma\cdot f_2(x^{\mu}) \leq B +\gamma\cdot(\opt(B)+\mu).
\end{align*}
\noindent
Since $\alg$ is an $\alpha$-approximation algorithm for $\Pi^{\WS}$, this implies that
\begin{align*}
  f_1(x^{\alg})+ \gamma\cdot f_2(x^{\alg}) \leq \alpha\cdot(f_1(x^{\mu}) + \gamma\cdot f_2(x^{\mu})) \leq \alpha\cdot(B +\gamma\cdot(\opt(B)+\mu)).
\end{align*}
As this holds for every $\mu>0$, we obtain that actually
\begin{align*}
	f_1(x^{\alg})+ \gamma\cdot f_2(x^{\alg}) \leq \alpha\cdot(B +\gamma\cdot\opt(B)).
\end{align*}
Moreover, as $\gamma=\bgamma=\frac{\epsilon\cdot B}{\opt(B)}\geq 0$ and $f_2(x^{\alg})\geq 0$, we can upper bound the first term on the left hand side by the right hand side to obtain that
\begin{align*}
  f_1(x^{\alg}) \leq \alpha\cdot(B +\bgamma\cdot\opt(B)) =  \alpha\cdot(B +\epsilon B) = \alpha\cdot(1+\epsilon)\cdot B.
\end{align*}
Similarly, using that $f_1(x^{\alg})\geq 0$ and $\gamma=\bgamma=\frac{\epsilon\cdot B}{\opt(B)}>0$ due to our assumptions on $\epsilon$, $B$, and $\opt(B)$, we obtain that
\begin{align*}
  f_2(x^{\alg}) & \leq \alpha\cdot\left(\frac{B}{\bgamma} + \opt(B)\right)=\alpha\cdot\left(\frac{\opt(B)}{\epsilon}+\opt(B)\right)\\
                & =\alpha\cdot\left(1+\frac{1}{\epsilon}\right)\cdot\opt(B),
\end{align*}
which proves the claim.
\end{proof}

Proposition~\ref{prop:approx} shows that, if we choose the ``correct'' weight~$\gamma=\bgamma=\frac{\epsilon\cdot B}{\opt(B)}$ for the weighted sum problem, then the solution returned by $\alg$
obtains a bicriteria approximation guarantee of $(\alpha\cdot(1+\epsilon),\alpha\cdot(1+\frac{1}{\epsilon}))$. Since we do not know the value~$\bgamma$ (as we do not know $\opt(B)$), we now
consider the approximation guarantee we can still obtain when using a weight slightly different from $\bgamma$:

\begin{prop}\label{prop:approx2}
Let $\bgamma=\frac{\epsilon\cdot B}{\opt(B)}$ as in Proposition~\ref{prop:approx}. If $\gamma\in[\frac{\bgamma}{1+\epsilon},(1+\epsilon)\bgamma]$, then the solution $x^{\alg}=x^{\alg}(\gamma)$
returned by $\alg$ for weight~$\gamma$ satisfies $f_1(x^{\alg})\leq \alpha\cdot(1+2\epsilon)\cdot B$ and $f_2(x^{\alg})\leq\alpha\cdot\left(1+\frac{2}{\epsilon}\right)\cdot\opt(B)$.
\end{prop}

\begin{proof}
As in the proof of Proposition~\ref{prop:approx}, the solution~$x^{\alg}$ satisfies
\begin{align*}
  f_1(x^{\alg})+ \gamma\cdot f_2(x^{\alg}) \leq \alpha\cdot(B +\gamma\cdot\opt(B)).
\end{align*}
Since $\gamma\leq (1+\epsilon)\bgamma=(1+\epsilon)\cdot\frac{\epsilon\cdot B}{\opt(B)}$ and the second summand $\gamma\cdot f_2(x^{\alg})$ is nonnegative, this implies that
\begin{align*}
  f_1(x^{\alg}) & \leq \underbrace{\alpha}_{\geq0}\cdot(B +\gamma\cdot\underbrace{\opt(B)}_{\geq0})\leq \alpha\cdot\left(B+(1+\epsilon)\cdot\frac{\epsilon\cdot B}{\opt(B)}\cdot\opt(B)\right) \\
                & = \alpha\cdot(B+\underbrace{(1+\epsilon)}_{\leq 2}\cdot\epsilon\cdot B) \leq \alpha\cdot(1+2\epsilon)\cdot B,
\end{align*}
where we have used that $1+\epsilon\leq 2$ due to the assumption that $\epsilon\leq 1$.

\medskip

Similarly, since $\gamma\geq \frac{\bgamma}{1+\epsilon}=\frac{\epsilon\cdot B}{(1+\epsilon)\cdot\opt(B)}>0$ and the first summand $f_1(x^{\alg})$ is nonnegative, we obtain
\begin{align*}
  f_2(x^{\alg}) & \leq \underbrace{\alpha}_{\geq0}\cdot\left(\frac{B}{\gamma} +\opt(B)\right)\leq \alpha\cdot\left(\frac{(1+\epsilon)\cdot\opt(B)\cdot B}{\epsilon\cdot B}+\opt(B)\right)\\
                & = \alpha\cdot\left(1+\frac{1+\epsilon}{\epsilon}\right)\cdot\opt(B) \leq \alpha\cdot\left(1+\frac{2}{\epsilon}\right)\cdot\opt(B),
\end{align*}
where we have again used that $1+\epsilon\leq 2$.
\end{proof}

Proposition~\ref{prop:approx2} shows that we achieve the desired approximation guarantee if we apply $\alg$ to the weighted sum problem for a weight $\gamma\in[\frac{\bgamma}{1+\epsilon},(1+\epsilon)\bgamma]$,
i.\,e., a weight that differs from $\bgamma$ by at most a factor $(1+\epsilon)$. Since we have $\LB(2)\leq \opt(B) \leq \UB(2)$ due to our assumptions on the objective functions, we have
$\bgamma=\frac{\epsilon\cdot B}{\opt(B)}\in[\frac{\epsilon\cdot B}{\UB(2)},\frac{\epsilon\cdot B}{\LB(2)}]$. Thus, we obtain the desired
$(\alpha\cdot(1+2\epsilon),\alpha\cdot(1+\frac{2}{\epsilon}))$-approximation if we apply $\alg$ to the values $(1+\epsilon)^i$ for all integers
$i=\left\lfloor\log_{(1+\epsilon)}\left(\frac{\epsilon\cdot B}{\UB(2)}\right)\right\rfloor$ up to $\left\lceil\log_{(1+\epsilon)}\left(\frac{\epsilon\cdot B}{\LB(2)}\right)\right\rceil$. This yields
a running time of
\begin{align*}
  \mathcal{O}\left(T_{\alg}\cdot\log_{(1+\epsilon)}\left(\frac{\UB(2)}{\LB(2)}\right)\right)=\mathcal{O}\left(T_{\alg}\cdot\frac{1}{\epsilon}\cdot\log_2\left(\frac{\UB(2)}{\LB(2)}\right)\right),
\end{align*}
where $T_{\alg}$ denotes the running time of $\alg$. This proves the following theorem:

\begin{thm}
Given $0<\epsilon\leq 1$ and a polynomial-time $\alpha$-approximation algorithm~$\alg$ for the weighted sum problem~$\Pi^{\WS}$, it is possible to obtain a bicriteria
$(\alpha\cdot(1+2\epsilon),\alpha\cdot(1+\frac{2}{\epsilon}))$-approximation algorithm for $\Pi^{\Budget}$ with running time
$\mathcal{O}\left(T_{\alg}\cdot\frac{1}{\epsilon}\cdot\log_2\left(\frac{\UB(2)}{\LB(2)}\right)\right)$, where $T_{\alg}$ denotes the running time of $\alg$.
\end{thm}

By setting $\epsilon\colonequals 1$, we obtain an algorithm yielding a fixed approximation guarantee in both objective functions:

\begin{cor}\label{cor:constapp}
	Given a polynomial-time $\alpha$-approximation algorithm~$\alg$ for the weighted sum problem~$\Pi^{\WS}$, it is possible to obtain a bicriteria
	$(3\alpha,3\alpha)$-approximation algorithm for $\Pi^{\Budget}$ with running time
	$\mathcal{O}\left(T_{\alg}\cdot\log_2\left(\frac{\UB(2)}{\LB(2)}\right)\right)$, where $T_{\alg}$ denotes the running time of $\alg$.
\end{cor}

\section{Improvements}\label{sec:running-time-imp}
In this section, we present two methods for improving our bicriteria approximation algorithm in the case that we are given an \emph{exact} polynomial-time algorithm~$\alg$ for the weighted sum problem (i.\,e., a $1$-approximation algorithm). In this case, we are able to obtain a bicriteria $(1+2\epsilon,1+\frac{2}{\epsilon})$-approximation algorithm with running time
$\mathcal{O}\left(T_{\alg}\cdot\log_2\left(\frac{1}{\epsilon}\cdot\log_2\left(\frac{\UB(2)}{\LB(2)}\right)\right)\right)$ by using binary search (whose combination with the usage of repeated powers of $(1+\epsilon)$ has already turned out to be useful within several other techniques for approximating bicriteria problems, cf.~\cite{papadimitriou2000approximability,vassilvitskii2005efficiently}). Moreover, if~$\alg$ additionally satisfies the necessary assumptions for using Megiddo's parametric search~\cite{Megiddo:rational-objectives,Megiddo:parallel}, we show in Subsection~\ref{subsec:parametric-search} that we can improve the approximation guarantee of the algorithm to $(1+\epsilon,1+\frac{1}{\epsilon})$ while at the same time achieving an improved running time.

\subsection{Binary Search}
Assume that we are given an \emph{exact} polynomial-time algorithm~$\alg$ for the weighted sum problem.
The idea for improving the running time is to perform a binary search on the integer values $i_{\min}\colonequals\left\lfloor\log_{(1+\epsilon)}\left(\frac{\epsilon\cdot B}{\UB(2)}\right)\right\rfloor$ up to
$i_{\max}\colonequals\left\lceil\log_{(1+\epsilon)}\left(\frac{\epsilon\cdot B}{\LB(2)}\right)\right\rceil$ instead of applying $\alg$ to each value in this range. The following lemma forms the basis for the
correctness of the binary search approach:

\begin{lem}\label{lem:binary-search}
For $i\in\{i_{\min},\dots,i_{\max}\}$, let $\gamma(i)\colonequals(1+\epsilon)^i$. If $i_{\min}\leq i < j \leq i_{\max}$, the solutions~$x^i\colonequals x^{\alg}(\gamma(i))$ and
$x^j\colonequals x^{\alg}(\gamma(j))$ returned by $\alg$ for the weights~$\gamma(i)$ and $\gamma(j)$, respectively, satisfy $f_1(x^i)\leq f_1(x^j)$ and $f_2(x^i)\geq f_2(x^j)$.
\end{lem}

\begin{proof}
Since we assume that $\alg$ is an exact algorithm for the weighted sum problem, $x^i$ is an optimal solution for the weighted sum problem with weight~$\gamma(i)$. Thus,
\begin{align}\label{eq:i-optimality}
  f_1(x^i) + \gamma(i)\cdot f_2(x^i) \leq f_1(x^j) + \gamma(i)\cdot f_2(x^j).
\end{align}
Similarly, since $x^j$ is an optimal solution for the weighted sum problem with weight~$\gamma(j)$, we have
\begin{align}\label{eq:j-optimality}
  f_1(x^j) + \gamma(j)\cdot f_2(x^j) \leq f_1(x^i) + \gamma(j)\cdot f_2(x^i).
\end{align}
Hence, we obtain that
\begin{align*}
  f_1(x^i) & \stackrel{\eqref{eq:i-optimality}}{\leq} f_1(x^j) + \gamma(i)\cdot f_2(x^j) - \gamma(i)\cdot f_2(x^i) \\
            & = f_1(x^j) + \gamma(j)\cdot f_2(x^j) + \bigl(\gamma(i)-\gamma(j)\bigr)\cdot f_2(x^j) - \gamma(i)\cdot f_2(x^i) \\
            & \stackrel{\eqref{eq:j-optimality}}{\leq} f_1(x^i) + \gamma(j)\cdot f_2(x^i) + \bigl(\gamma(i)-\gamma(j)\bigr)\cdot f_2(x^j) - \gamma(i)\cdot f_2(x^i) \\
            & = f_1(x^i) + \bigl(\gamma(j)-\gamma(i)\bigr)\cdot\bigl(f_2(x^i)-f_2(x^j)\bigr).
\end{align*}
Subtracting $f_1(x^i)$ on both sides and using that $(\gamma(j)-\gamma(i))=(1+\epsilon)^j - (1+\epsilon)^i>0$ since $j>i$ shows that $f_2(x^i)\geq f_2(x^j)$, which proves the claimed relation of the $f_2$-costs.

\medskip
\noindent
To show the claimed relation of the $f_1$-cost, we use that
\begin{align*}
  f_1(x^i) + \gamma(i)\cdot f_2(x^i) & \stackrel{\eqref{eq:i-optimality}}{\leq} f_1(x^j) + \underbrace{\gamma(i)}_{\geq0}\cdot \underbrace{f_2(x^j)}_{\leq f_2(x^i)} \leq f_1(x^j) + \gamma(i)\cdot f_2(x^i),
\end{align*}
so subtracting $\gamma(i)\cdot f_2(x^i)$ from both sides shows that $f_1(x^i)\leq f_1(x^j)$ as claimed.
\end{proof}

Lemma~\ref{lem:binary-search} shows that we can perform a binary search on the values in $\{i_{\min},\dots,i_{\max}\}$ as follows:
\begin{description}
\item If we run $\alg$ with weight~$\gamma(i)$ for some $i\in\{i_{\min},\dots,i_{\max}\}$ and the resulting solution~$x^i$ satisfies $f_1(x^i)>(1+2\epsilon)\cdot B$, then Lemma~\ref{lem:binary-search} shows
      that $f_1(x^j)\geq f_1(x^i)>(1+2\epsilon)\cdot B$ for all $j>i$. Hence, no value~$j>i$ yields a solution that obtains the desired approximation guarantee in the first objective function and we only need
      to consider the values smaller than~$i$.
\item If, on the other hand, the resulting solution~$x^i$ satisfies $f_1(x^i)\leq (1+2\epsilon)\cdot B$, then Lemma~\ref{lem:binary-search} (with the roles of $i$ and $j$ exchanged) shows that
      $f_2(x^j)\geq f_2(x^i)$ for all $j<i$. Hence, the solution~$x^i$ obtains the desired approximation guarantee in the first objective function and no value~$j<i$ can yield a solution with a better value
      in the second objective function. Consequently, we only need to consider the values larger than~$i$.
\end{description}

When using binary search as described above on the values in $\{i_{\min},\dots,i_{\max}\}$, we only need $\mathcal{O}\left(\log_2\left(\log_{(1+\epsilon)}\left(\frac{\UB(2)}{\LB(2)}\right)\right)\right)$ calls
of $\alg$, which leads to a total running time of
\begin{align*}
  \mathcal{O}\left(T_{\alg}\cdot\log_2\left(\log_{(1+\epsilon)}\left(\frac{\UB(2)}{\LB(2)}\right)\right)\right)
  =\mathcal{O}\left(T_{\alg}\cdot\log_2\left(\frac{1}{\epsilon}\cdot\log_2\left(\frac{\UB(2)}{\LB(2)}\right)\right)\right),
\end{align*}
where $T_{\alg}$ denotes the running time of $\alg$. This proves the following theorem:

\begin{thm}
Given $0<\epsilon\leq 1$ and a polynomial-time exact algorithm~$\alg$ for the weighted sum problem~$\Pi^{\WS}$, it is possible to obtain a bicriteria\\
$(1+2\epsilon,1+\frac{2}{\epsilon})$-approximation algorithm for $\Pi^{\Budget}$ with running time\\
$\mathcal{O}\left(T_{\alg}\cdot\log_2\left(\frac{1}{\epsilon}\cdot\log_2\left(\frac{\UB(2)}{\LB(2)}\right)\right)\right)$, where $T_{\alg}$ denotes the running time of~$\alg$.
\end{thm}

Note that, if we are only given an $\alpha$-approximation algorithm for the weighted sum problem with $\alpha>1$, we cannot perform a binary search as above. The reason is that the solutions returned by the approximation algorithm need not fulfill inequalities~\eqref{eq:i-optimality} and \eqref{eq:j-optimality}, so Lemma~\ref{lem:binary-search} does not hold in this more general case.

\subsection{Parametric Search}\label{subsec:parametric-search}
We now show how our algorithm can be improved to yield a better approximation guarantee of $(1+\epsilon,1+\frac{1}{\epsilon})$ as well as an improved running time by using Megiddo's parametric search technique~\cite{Megiddo:rational-objectives,Megiddo:parallel} if the polynomial-time exact algorithm~$\alg$ for the weighted sum problem~$\Pi^{\WS}$ satisfies certain additional assumptions. Specifically, we need to assume that, when $\alg$ is applied to the weighted sum problem~$\Pi^{\WS}(\gamma)$, the value~$\gamma$ is only involved in additions, comparisons, and multiplications by constants during the execution of $\alg$. Additionally, we assume that the solution returned by $\alg$ depends on the value of~$\gamma$ only via the outcomes of the comparisons involving~$\gamma$ that are made during the execution of $\alg$, so the returned solution is independent of~$\gamma$ as long as the outcomes of these comparisons do not change.

\medskip

Recall that, by Proposition~\ref{prop:approx}, $\alg$ returns a solution obtaining an approximation guarantee of $(1+\epsilon,1+\frac{1}{\epsilon})$ if we choose the ``correct'' weight~$\gamma=\bgamma$. The idea is now to search for a value~$\hat{\gamma}$ (which need not necessarily be equal to~$\bgamma$) for which the corresponding solution returned by~$\alg$ yields this approximation guarantee parametrically as described by Megiddo's general method: We run a master copy of the algorithm~$\alg$ with weight~$\gamma$ (i.\,e., objective function $f_1+\gamma\cdot f_2$) and keep the value~$\gamma$ as a parameter. The execution of this algorithm will proceed in the same way for all possible values of~$\gamma$ except for points in time where a comparison is made that involves~$\gamma$. In this case, the outcome of the comparison (and, thus, the path of computation taken by the algorithm) may depend on $\gamma$.
However, in the input data given to $\alg$, the value~$\gamma$ appears only in the objective function $f_1+\gamma\cdot f_2$, whose value at any point depends linearly on~$\gamma$, and the only operations (besides comparisons) performed during the execution of~$\alg$ that involve~$\gamma$ are additions and multiplications by constants, which maintain linearity in~$\gamma$. Hence, any comparison during the execution of $\alg$ that involves~$\gamma$ will be between two \emph{linear} functions of~$\gamma$. Thus, there will be at most one critical value (say~$\gamma'$) such that, for $\gamma\leq \gamma'$, one of the functions is greater than or equal to the other and, for $\gamma\geq\gamma'$, this relation is reversed. Hence, the comparison may partition the real line into two subintervals to the left and to the right of~$\gamma'$ such that the outcome of the comparison is identical on each subinterval.

\medskip

The central observation is now that, using the monotonicity of the values of~$f_1$ and $f_2$ when applied to the optimal solutions of the weighted sum problem for different values of~$\gamma$ shown in Lemma~\ref{lem:binary-search}, we can resolve the comparison (i.\,e., decide for either the subinterval to the left or the right of the critical value~$\gamma'$) by one application of~$\alg$ with weight~$\gamma'$: If~$x'$ denotes the solution returned by $\alg$ for $\gamma'$ and $f_1(x')>(1+\epsilon)\cdot B$, then the monotonicity observed in Lemma~\ref{lem:binary-search} shows that any solution~$x$ obtained from~$\alg$ for a value $\gamma>\gamma'$ satisfies $f_1(x)\geq f_1(x')>(1+\epsilon)\cdot B$. Thus, no value~$\gamma>\gamma'$ yields a solution that obtains the desired approximation guarantee in the first objective function and we can conclude that $\bgamma\leq \gamma'$. In particular, we know that the subinterval of the real line to the left of $\gamma'$ still contains a value~$\hat{\gamma}$ as desired. Thus, we can choose this subinterval and continue with the execution of the master copy of~$\alg$ by resolving the comparison as in this subinterval.

Similarly, if $f_1(x')\leq (1+\epsilon)\cdot B$, then the monotonicity observed in Lemma~\ref{lem:binary-search} shows that any solution~$x$ obtained from~$\alg$ for a value $\gamma<\gamma'$ satisfies $f_2(x)\geq f_2(x')$. Hence, $x'$ obtains the desired approximation guarantee in the first objective function and no value~$\gamma<\gamma'$ yields a solution with a better value in the second objective function. Thus, we can conclude that there exists some value~$\hat{\gamma}$ with $\hat{\gamma}\geq \gamma'$ that yields a solution obtaining the desired approximation guarantees in both objective functions, so we can continue the execution of the master copy of~$\alg$ by resolving the comparison as in the subinterval to the right of~$\gamma'$ (note that, in this case, we cannot be sure that $\bgamma\geq \gamma'$, so the value~$\hat{\gamma}$ we obtain at the end might be different from~$\bgamma$).

\medskip

When the execution of the master algorithm ends, we obtain an interval~$[a,b]$ and a solution~$x$ of the problem instance such that~$\alg$ returns the solution~$x$ for all $\gamma\in[a,b]$ (since all comparisons involving~$\gamma$ have the same outcome for all~$\gamma\in[a,b]$, the solution is independent of the choice of~$\gamma\in[a,b]$ by assumption). Moreover, an inductive application of the above arguments shows that $[a,b]$ still contains a value~$\hat{\gamma}$ whose corresponding solution of the weighted sum problem returned by~$\alg$ (which must be~$x$) yields the desired approximation guarantee of $(1+\epsilon,1+\frac{1}{\epsilon})$. Consequently, the solution~$x$ is as desired.

\medskip

Since we need to call~$\alg$ once for every comparison made during the execution of the master copy of~$\alg$, the running time of the procedure is determined by the number of comparisons in the master copy of~$\alg$ multiplied with the running time of~$\alg$. This proves the following theorem:

\begin{thm}\label{thm:meggido}
Suppose that we are given $0<\epsilon\leq 1$ and a polynomial-time exact algorithm~$\alg$ for the weighted sum problem~$\Pi^{\WS}$ such that, when $\alg$ is applied to the weighted sum problem~$\Pi^{\WS}(\gamma)$, the value~$\gamma$ is only involved in additions, comparisons, and multiplications by constants and the solution returned by $\alg$ depends on the value of~$\gamma$ only via the outcomes of the comparisons involving~$\gamma$ that are made during the execution of the algorithm. Then it is possible to obtain a bicriteria
$(1+\epsilon,1+\frac{1}{\epsilon})$-approximation algorithm for~$\Pi^{\Budget}$ with running time~$\mathcal{O}\left(c\cdot T_{\alg}\right)$, where $T_{\alg}$ denotes the running time of~$\alg$ and~$c$ denotes the number of comparisons performed during one execution of~$\alg$.
\end{thm}

We remark that, in the case that~$\alg$ is a \emph{strongly} polynomial-time algorithm, the running time we obtain will also be strongly polynomial (it is at most $\mathcal{O}\left((T_{\alg})^2\right)$). We also remark that, in many cases, it is possible to improve the running time further by exploiting parallelism within~$\alg$ and resolving several independent comparisons at once by using binary search on the corresponding critical values. We refer to \cite{Megiddo:rational-objectives,Megiddo:parallel} for the details.

\section{Obtaining an Approximate Pareto Curve}\label{sec:Pareto-curve}

We now show how the method used for obtaining our bicriteria approximation algorithm for~$\Pi^{\Budget}$ can be adjusted in order to compute an approximate Pareto curve for~$\Pi$. As in Section~\ref{sec:bicriteria-approx}, we again assume
that we are given $0<\epsilon\leq1$ and a polynomial-time $\alpha$-approximation algorithm~$\alg$ for the weighted sum problem~$\Pi^{\WS}$.

\medskip

Suppose that~$\tilde{x}$ is any feasible solution and let $\tilde{B}\colonequals f_1(\tilde{x})$, $\bgamma(\tilde{B})\colonequals\frac{\epsilon\cdot \tilde{B}}{\opt(\tilde{B})}$, where, as before,
$\opt(\tilde{B})=\inf\{f_2(x):x\in S, f_1(x)\leq \tilde{B}\}$ denotes the minimum $f_2$-cost among the feasible solutions with $f_1$-cost at most~$\tilde{B}$. By Proposition~\ref{prop:approx2}, we obtain a solution~$x$ with
$f_1(x)\leq \alpha\cdot(1+2\epsilon)\cdot f_1(\tilde{x})$ and $f_2(x)\leq \alpha\cdot(1+\frac{2}{\epsilon})\cdot f_2(\tilde{x})$ if we apply $\alg$ to the weighted sum problem for some weight
$\gamma\in\left[\frac{\bgamma(\tilde{B})}{1+\epsilon},(1+\epsilon)\bgamma(\tilde{B})\right]$. Moreover, since $\LB(2)\leq \opt(\tilde{B})\leq UB(2)$, we have
$\bgamma(\tilde{B})\in\left[\frac{\epsilon\cdot\tilde{B}}{\UB(2)},\frac{\epsilon\cdot\tilde{B}}{\LB(2)}\right]$. As $\LB(1)\leq \tilde{B}=f_1(\tilde{x})\leq UB(1)$ for all feasible solutions~$\tilde{x}$,
this shows that, by applying $\alg$ to the weighted sum problem with weight $\gamma=(1+\epsilon)^i$ for $i=\left\lfloor\log_{(1+\epsilon)}\left(\frac{\epsilon\cdot \LB(1)}{\UB(2)}\right)\right\rfloor$ up to $\left\lceil\log_{(1+\epsilon)}\left(\frac{\epsilon\cdot \UB(1)}{\LB(2)}\right)\right\rceil$ and forming the union of all returned solutions, we obtain an $(\alpha\cdot(1+2\epsilon),\alpha\cdot(1+\frac{2}{\epsilon}))$-approximate Pareto curve. The total running time obtained is
\begin{align*}
  \mathcal{O}\left(T_{\alg}\cdot\log_{(1+\epsilon)}\left(\frac{\UB(1)\cdot\UB(2)}{\LB(1)\cdot\LB(2)}\right)\right)
  =\mathcal{O}\left(T_{\alg}\cdot\frac{1}{\epsilon}\cdot\log_2\left(\frac{\UB(1)\cdot\UB(2)}{\LB(1)\cdot\LB(2)}\right)\right),
\end{align*}
where $T_{\alg}$ denotes the running time of $\alg$. This proves the following theorem:

\begin{thm}\label{thm:approx-Pareto-curve}
Given $0<\epsilon\leq 1$ and a polynomial-time $\alpha$-approximation algorithm~$\alg$ for the weighted sum problem~$\Pi^{\WS}$, it is possible to compute an
$(\alpha\cdot(1+2\epsilon),\alpha\cdot(1+\frac{2}{\epsilon}))$-approximate Pareto curve in time\\
$\mathcal{O}\left(T_{\alg}\cdot\frac{1}{\epsilon}\cdot\log_2\left(\frac{\UB(1)\cdot\UB(2)}{\LB(1)\cdot\LB(2)}\right)\right)$, where $T_{\alg}$ denotes the running time of $\alg$.
\end{thm}

Analogously to Corollary~\ref{cor:constapp}, setting $\epsilon\colonequals 1$ shows that we can obtain a $(3\alpha,3\alpha)$-approximate Pareto curve in polynomial time whenever a polynomial-time $\alpha$-approximation algorithm~$\alg$ for the weighted sum problem is available.

\medskip

We remark that we can obtain a better approximation to the Pareto curve in case that there exists a polynomial-time \emph{parametric} ($\alpha$-approximation) algorithm for the weighted sum problem~$\Pi^{\WS}$. Such an algorithm computes a sequence of ($\alpha$-approximate) solutions for the weighted sum problem~$\Pi^{\WS}(\gamma)$ for all $\gamma>0$ simultaneously in polynomial time (the output usually consists of a sequence of polynomially many intervals for~$\gamma$ and a corresponding ($\alpha$-approximate) solution for each interval). Since Proposition~\ref{prop:approx} shows that, for any feasible solution~$\tilde{x}$, an $\alpha$-approximate solution~$x$ of $\Pi^{\WS}(\bgamma(\tilde{B}))$ satisfies $f_1(x)\leq \alpha\cdot(1+\epsilon)\cdot f_1(\tilde{x})$ and $f_2(x)\leq \alpha\cdot(1+\frac{1}{\epsilon})\cdot f_2(\tilde{x})$, it follows immediately that the set of solutions returned by a parametric $\alpha$-approximation algorithm is an $(\alpha\cdot(1+\epsilon),\alpha\cdot(1+\frac{1}{\epsilon}))$-approximate Pareto curve. The running time for obtaining this set is equal to the running time of the parametric algorithm.

\section{Applications}\label{sec:applications}
Our results are applicable to a vast variety of minimization problems, especially combinatorial problems. In this section, we provide some examples of such problems and examine the specific results obtained from our general method. Here, we focus on the results concerning the computation of approximate Pareto curves. The running time of our method for the corresponding budget-constrained problems can easily be calculated from our general theorems by using the stated running times of the approximation algorithms used for the weighted sum problems.

\medskip

The presented method applies to all bicriteria minimization problems with positive-valued, polynomially computable objective functions for which an $\alpha$-approximation algorithm for the weighted sum problem exists. For problems in which the two objective functions are of the same type, such an approximation for the weighted sum problem can usually be obtained by using an approximation algorithm for the corresponding unicriterion problem. In many graph problems, e.g., the vertices or arcs of the graph are weighted by positive rational numbers and the feasible solutions correspond to subsets of the vertices or arcs, whose objective function values are given by the sums of the corresponding weights of these vertices/arcs. The weighted sum problem then simply corresponds to the unicriterion problem in which the single weight of each vertex/arc is given by the weighted sum of its two weights from the bicriteria problem.

\medskip

The assumptions needed to apply our method are rather weak and the method can be employed on a large class of bicriteria minimization problems. In particular, our method is applicable to the bicriteria versions of \emph{$\textsf{APX}$-complete} problems, which do not admit multicriteria polynomial-time approximation schemes unless $\textsf{P}=\textsf{NP}$ since the corresponding gap problems cannot be solved in polynomial time.\footnote{An $\textsf{NP}$-optimization problem is in the class~$\textsf{APX}$ if it admits a polynomial-time constant-factor approximation algorithm, and it is called \emph{$\textsf{APX}$-complete} if it is contained in $\textsf{APX}$ and every problem in $\textsf{APX}$ can be reduced to it via a PTAS reduction. It is known that $\textsf{APX}$-complete problems do not admit a PTAS unless $\textsf{P}=\textsf{NP}$.} Hence, the results of Papadimitriou and Yannakakis~\cite{papadimitriou2000approximability} and the succeeding papers~\cite{vassilvitskii2005efficiently,diakonikolas2009small,Bazgan+etal:min-pareto} cannot be used for these problems. Thus, the only general method known for computing approximate Pareto curves for these problems is by combining the results of Glaser et al.~\cite{glasser2010approximability,Glasser+etal:multi-hardness} as mentioned in the introduction, for which no running time analysis is available.

\medskip

Table~\ref{table:results} provides an exemplary list of problems for which our method can be used in order to compute approximate Pareto curves and compares the obtained running times and approximation guarantees to the previously best known methods for computing approximate Pareto curves for these problems. For precise definitions of the specific problems, we refer to the references provided in the table. For most of the mentioned problems, our method provides the first approximation result with an explicit analysis of the running time while at the same time obtaining approximation guarantees that are comparable to the best previously known methods. For other problems, for which an MFPTAS exists, our method yields larger approximation ratios but significantly faster running times than the fastest known MFPTAS.

\begin{sidewaystable}
\centering
\setlength{\tabcolsep}{1pt}
\begin{tabular}{|C{2.5cm}|C{3.5cm}|C{3cm}|C{3.5cm}|C{4cm}|C{5.9cm}|}
\hline
\bf{Problem} & \bf{Prev. Approx. Guarantees} & \bf{Running Time} & \bf{Our Approx. Guarantee} & \bf{Alg. for $\Pi^{\WS}$} & \bf{Our Running Time} \\
\hline
\multicolumn{6}{|c|}{\textbf{$\textsf{APX}$-complete problems}} \\
\hline
Bicrit. metric TSP &
$(2,2)$~\cite{glasser2009improved2, glasser2009improved} & not specified & $(\frac{3}{2}+3\epsilon, \frac{3}{2}+\frac{3}{\epsilon})$ [Thm.~\ref{thm:approx-Pareto-curve}] & $\frac{3}{2}$-approxi\-ma\-tion~\cite{Christofides:TSP} with running time $\mathcal{O}(n^3)$ & $\mathcal{O}\left(n^3\cdot\frac{1}{\epsilon}\cdot\log_2\left(\frac{\UB(1)\cdot\UB(2)}{\LB(1)\cdot\LB(2)}\right)\right)$ \\
\hline
Bicrit. min. weight vertex cover & $(4+4\epsilon,4+4\epsilon)$ using~\cite{glasser2010approximability,Glasser+etal:multi-hardness} & not specified & $(2+4\epsilon, 2+\frac{4}{\epsilon})$ [Thm.~\ref{thm:approx-Pareto-curve}] & $2$-approxi\-ma\-tion~\cite{bar1981linear} with running time $\mathcal{O}(m)$ & $\mathcal{O}\left(m\cdot\frac{1}{\epsilon}\cdot\log_2\left(\frac{\UB(1)\cdot\UB(2)}{\LB(1)\cdot\LB(2)}\right)\right)$ \\
\hline
Bicrit. min. $k$-spanning tree & $(4+4\epsilon,4+4\epsilon)$ using~\cite{glasser2010approximability,Glasser+etal:multi-hardness} & not specified & $(2+4\epsilon, 2+\frac{4}{\epsilon})$ [Thm.~\ref{thm:approx-Pareto-curve}] & $2$-approxi\-ma\-tion~\cite{garg2005saving} with running time $\mathcal{O}(m\cdot n^4\cdot\log_2 n)$ & $\mathcal{O}\Bigl(m\cdot n^4\cdot\log_2 n\cdot\frac{1}{\epsilon}\cdot$\newline $\left. \log_2\left(\frac{\UB(1)\cdot\UB(2)}{\LB(1)\cdot\LB(2)}\right)\right)$ \\
\hline
Bicrit. min. weight edge dom. set & $(4+4\epsilon,4+4\epsilon)$ using~\cite{glasser2010approximability,Glasser+etal:multi-hardness} & not specified & $(2+4\epsilon, 2+\frac{4}{\epsilon})$ [Thm.~\ref{thm:approx-Pareto-curve}] & $2$-approxi\-ma\-tion~\cite{fujito2002} (running time not specified) & $\mathcal{O}\left(T_{\alg}\cdot\frac{1}{\epsilon}\cdot\log_2\left(\frac{\UB(1)\cdot\UB(2)}{\LB(1)\cdot\LB(2)}\right)\right)$, where $T_{\alg}$ is the running time of the alg. from~\cite{fujito2002}. \\
\hline
Bicrit. min. weight Steiner tree & $((\ln(4)+\delta)(2+2\epsilon),$\newline $(\ln(4)+\delta)(2+2\epsilon))$ for any $\delta>0$ using~\cite{glasser2010approximability,Glasser+etal:multi-hardness} & not specified & $((\ln(4)+\delta)(1+2\epsilon),$\newline $(\ln(4)+\delta)(1+\frac{2}{\epsilon}))$ for any $\delta>0$ [Thm.~\ref{thm:approx-Pareto-curve}] & $(\ln(4)+\delta)$-appro\-xi\-ma\-tion~\cite{byrka2010improved} (running time not specified) & $\mathcal{O}\left(T_{\alg}\cdot\frac{1}{\epsilon}\cdot\log_2\left(\frac{\UB(1)\cdot\UB(2)}{\LB(1)\cdot\LB(2)}\right)\right)$, where $T_{\alg}$ is the running time of the alg. from~\cite{byrka2010improved}. \\
\hline
\multicolumn{6}{|c|}{\textbf{Other problems}} \\
\hline
Bicrit. min. $s$-$t$-cut & $(2+\epsilon,2+\epsilon)$ using~\cite{glasser2010approximability,Glasser+etal:multi-hardness} and no MFPTAS unless $\textsf{P}=\textsf{NP}$~\cite{papadimitriou2000approximability} & not specified & $(1+2\epsilon,1+\frac{2}{\epsilon})$ [Thm.~\ref{thm:approx-Pareto-curve}] & max. flow algs. from~\cite{Orlin:STOC13,King+etal:max-flow} with running time $\mathcal{O}(nm)$ & $\mathcal{O}\left(n\cdot m\cdot\frac{1}{\epsilon}\cdot\log_2\left(\frac{\UB(1)\cdot\UB(2)}{\LB(1)\cdot\LB(2)}\right)\right)$ \\
\hline
Bicrit. shortest $s$-$t$-path & MFPTAS~\cite{Tsaggouris+Zaroliagis:mult-shortest-path} & $\mathcal{O}(n^2\cdot m\cdot\frac{1}{\epsilon}\cdot\log_2(n\cdot C^{\max}))$ & $(1+\epsilon,1+\frac{1}{\epsilon})$ [Thm.~\ref{thm:approx-Pareto-curve}] and use of parametric alg. for solving~$\Pi^{\WS}$ & parametric shortest path alg. from~\cite{Young+etal:parametric-sp} with running time $\mathcal{O}(m\cdot n + n^2\cdot\log_2 n)$ & $\mathcal{O}(m\cdot n + n^2\log_2 n)$ \newline (strongly polynomial, equal to parametric alg. used) \\
\hline
Bicrit. min spanning tree & MFPTAS~\cite{vassilvitskii2005efficiently,diakonikolas2009small,Bazgan+etal:min-pareto} & $\mathcal{O}(\left(\frac{1}{\epsilon}\right)^{\frac{1}{\epsilon}}\cdot n^3\cdot |P^*_{\epsilon}|)$ or\newline $\mathcal{O}\left(m\cdot n^5\cdot|P^*_{\epsilon}|\cdot\right.$\newline $\left.\tau(\small\lfloor\frac{n-1}{\epsilon}\rfloor,\lfloor\frac{n-1}{\epsilon}\rfloor)\right)$ & $(1+\epsilon,1+\frac{1}{\epsilon})$ [Thm.~\ref{thm:approx-Pareto-curve}] and use of parametric alg. for solving~$\Pi^{\WS}$ & parametric min. spanning tree alg. from~\cite{Fernandez-Baca+etal:SWAT96} that runs in $\mathcal{O}(m\cdot n \cdot \log_2 n)$ time & $\mathcal{O}(m\cdot n \cdot \log_2 n)$ \newline (strongly polynomial, equal to parametric alg. used) \\
\hline
\end{tabular}
\caption{Comparison of our results on the computation of approximate Pareto curves to previous results. $|P^*_{\epsilon}|$ denotes the cardinality of the smallest $\epsilon$-Pareto curve and $\tau(a,b)$ denotes the time to multiply two polynomials of maximum degrees~$a$ and~$b$.}\label{table:results}
\end{sidewaystable}

\section{Maximization Problems}\label{sec:maximization}
An obvious question is whether our method can be adapted for bicriteria maximization problems with positive-valued, polynomially computable objective functions.\footnote{Note that, due to the assumption that the values of the objective functions must be positive, our method cannot be applied to such problems by simply reversing the sign of the objective functions and replacing the maximization by a minimization of all objective functions.} In this section, we show that a bicriteria approximation algorithm providing the same approximation guarantees as our method cannot exist for the budget-constrained version of general bicriteria maximization problems unless $\textsf{P}=\textsf{NP}$.\footnote{Note that, in the case of a maximization problem, the budget constraint imposes a \emph{lower} bound on the value of the first objective function instead of an upper bound. We still refer to this constraint as a budget constraint here in order keep the terminology consistent with the one used for minimization problems.}
To do so, we use a slight modification of a construction utilized by Gla{\ss}er et. al~\cite{Glasser+etal:multi-hardness}. They considered a restricted version of the bicriteria maximum weight clique problem to
show that certain translations of approximability results from the weighted sum version of an optimization problem to the multicriteria version that work for minimization problems are not possible in
general for maximization problems.

\medskip

The bicriteria maximum weight clique problem ($\twoCLIQUE$) is defined as follows:

\begin{prob}[$\twoCLIQUE$]
	Given an undirected graph~$G=(V,E)$ and a pair of positive rational weights~$(w_1(v),w_2(v))\in\mathbb{Q}_{>0}^2$ for each vertex~$v\in V$, compute a clique (i.\,e., a subset~$V'\subseteq V$ of the vertex
	set such that every two vertices in $V'$ are connected by an edge) maximizing $\sum_{v\in V'} (w_1(v),w_2(v))$ among all cliques in~$G$.
\end{prob}

Similar to \cite{Glasser+etal:multi-hardness}, we consider a restricted version of $\twoCLIQUE$ (which we denote by $\twoCLIQUE_{\text{rest}}$) in which an instance consists of an arbitrary graph~$G$ with
weights~$(1,1)$ for all vertices and two additional isolated vertices~$x,y$ that have weights~$(2n+1,\frac{1}{n})$ and $(\frac{1}{n},2n+1)$, respectively, where $n=|V|$.

Note that (when excluding the empty clique that yields objective value zero in both objective functions from the set of feasible solutions) the problem $\twoCLIQUE_{\text{rest}}$ satisfies all our assumptions
(we can choose $\LB(i)\colonequals \frac{1}{n}$ and $\UB(i)\colonequals 2n+1$ for $i=1,2$) except that the objective functions are to be maximized instead of minimized. Moreover, the weighted sum problem
corresponding to $\twoCLIQUE_{\text{rest}}$ can be solved optimally in polynomial time for any $\gamma>0$ as follows: The algorithm simply outputs $\{x\}$ if $\gamma\leq 1$ and $\{y\}$ otherwise.
We now show, however, that no polynomial-time $(\alpha,\beta)$-approximation algorithm can exist for the budget-constrained version of $\twoCLIQUE_{\text{rest}}$ unless $\textsf{P}=\textsf{NP}$:

\begin{thm}
The budget-constrained version of $\twoCLIQUE_{\text{rest}}$ cannot be approximated in polynomial time with an approximation guarantee~$(\alpha,\beta)$ for any two constants~$\alpha,\beta\geq 1$ unless $\textsf{P}=\textsf{NP}$.
\end{thm}

\begin{proof}
We show that the existence of a polynomial-time $(\alpha,\beta)$-approximation algorithm~$\alg$ for the budget-constrained version of $\twoCLIQUE_{\text{rest}}$ implies the existence of a $\beta$-approximation algorithm for the (unweighted) unicriterion maximum clique problem (which we denote by $\CLIQUE$ in the following). Since approximating $\CLIQUE$ within any constant factor is $\textsf{NP}$-hard, this will show the claim.

So assume that an $(\alpha,\beta)$-approximation algorithm~$\alg$ for the budget-constrained version of $\twoCLIQUE_{\text{rest}}$ exists and consider an arbitrary instance of $\CLIQUE$ consisting of an undirected graph $G=(V,E)$ with $n$~vertices.
We consider the instance of the budget-constrained version of $\twoCLIQUE_{\text{rest}}$ defined by the graph~$G$ (with weights~$(1,1)$ for all vertices) together with the two isolated vertices~$x,y$ with weight~$(2n+1,\frac{1}{n})$ and $(\frac{1}{n},2n+1)$, respectively, and a lower bound of~$B\colonequals 1$ on the value of the first objective function. Without loss of generality, we assume that $\frac{1}{n}<
\min\{\frac{1}{\alpha},\frac{k}{\beta}\}$, where~$k$ denotes the maximum size of a clique within $G$ (this can be achieved by adding additional
isolated vertices with weights~$(1,1)$ to $G$, which increases~$n$, but does not change~$k$). Let~$W$ denote the clique returned by $\alg$ when applied to this instance of the budget-constrained problem. Since all vertices in $V$ have weight~$(1,1)$, we know that a feasible solution with objective function value~$k$ exists for the instance. Hence, the clique~$W$ returned by $\alg$ must have total weight at least $\frac{1}{\alpha}$ with respect to the first objective function and at least $\frac{k}{\beta}$ with respect to the second objective function. In particular, since $y$ has weight~$\frac{1}{n}<\frac{1}{\alpha}$ in the first objective function and $x$ has weight~$\frac{1}{n}<\frac{k}{\beta}$ in the second objective function, we have and $W\neq \{y\}$ and $W\neq \{x\}$. Hence, we must have $W\subseteq V$ and, since all vertices in $V$ have weight~$1$ in the second objective function, we obtain that $W$ is a clique of size at least $\frac{k}{\beta}$ in $G$.
\end{proof}

\section{Conclusion}\label{sec:conclusion} 
This paper provides a novel algorithm for approximating bicriteria minimization problems. It relies on the exact or approximate solution of the weighted sum problem corresponding to the bicriteria problem. Interestingly, in two variants, it is either possible to approximate the budget-constrained problem (also known as the $\epsilon$-constraint scalarization) or the whole Pareto curve. Our algorithm extends the state of the art in several ways: (a) it is generally applicable to a large variety of bicriteria minimization problems, (b) it illuminates the relationship between bicriteria and scalarized unicriterion minimization problems, (c) it provides first ever approximations with an explicit running time analysis for some minimization problems or (d), for other problems, it provides alternative approximations or approximations with an improved running time when compared to the state of the art, and (e) it highlights some principle differences between bicriteria maximization and minimization problems.

\medskip

The work at hand focuses exclusively on bicriteria minimization problems. Three directions of research immediately suggest themselves: First, the question about generalizations to more than two objective functions is certainly interesting. Second, the interplay of other scalarizations of bicriteria problems with respect to approximability should be deeper studied in the future. Third, in view of our result for maximization problems, the possibility of approximating general or a particular class of maximization problems is still an open research question.

\section{Acknowledgement}
The research of Pascal Halffmann was supported by DFG grant RU 1524/4-1, Stefan Ruzika acknowledges support by DAAD project 57128839, and David Willems' activities were funded by BMBF project 13N12825.

\bibliographystyle{elsarticle-num}
\bibliography{Literatur}

\section*{Appendix~$1$}
We show that we can relax the assumption of strictly positive objective values by only assuming nonnegativity of the objective values in the case that we can compute positive rational numbers
$\LB(I,i),\UB(I,i)$ for $i=1,2$ such that $\LB(I,i)\leq f_i(I,x) \leq \UB(I,i)$ for every feasible solution~$x\in S(I)$ with $f_i(I,x)>0$ in polynomial time.

\medskip

To this end, note that the image of the Pareto curve for a given instance~$I$ of $\Pi$ can contain at most one point~$(a,0)$ with second component zero and at most one point~$(0,b)$ with first component zero.
We show that, by applying the $\alpha$-approximation algorithm~$\alg$ to the corresponding instance of the weighted sum problem~$\Pi^{\WS}(\gamma)$ for two suitably chosen values of $\gamma$, we can
compute solutions whose images approximate these points~$(a,0)$ and $(0,b)$ within the desired approximation guarantees in polynomial time. 
For the rest of the Pareto curve, we can then apply our method using the assumption that the objective values are strictly positive.

\medskip

So first assume that a point with image~$(a,0)$ exists on the Pareto curve for a given instance~$I$. We apply $\alg$ to the corresponding instance of the weighted sum
problem~$\Pi^{\WS}(\gamma)$ for some $\gamma>\alpha\cdot\frac{\UB(I,1)}{\LB(I,2)}$. Since $\alg$ is an $\alpha$-approximation algorithm for the weighted sum problem, the resulting
solution~$x^{\alg}=x^{\alg}(\gamma)$ satisfies
\begin{align*}
  f_1(I,x^{\alg})+ \gamma\cdot f_2(I,x^{\alg}) \leq \alpha\cdot(a +\gamma\cdot 0) = \alpha\cdot a. 
\end{align*}
As $\gamma\cdot f_2(I,x^{\alg})\geq 0$, this implies that $f_1(I,x^{\alg})\leq \alpha\cdot a\leq \alpha\cdot(1+2\epsilon)\cdot a$. Hence, $f(I,x^{\alg})$ approximates $(a,0)$ with the desired approximation guarantee in the first component. Moreover, since also $f_1(I,x^{\alg})\geq 0$, we also obtain that $\gamma\cdot f_2(I,x^{\alg}) \leq \alpha\cdot a$, so by our choice of $\gamma$ and since $a\leq \UB(I,1)$:
\begin{align*}
  f_2(I,x^{\alg}) \leq \frac{\alpha\cdot a}{\gamma} < \frac{a \cdot \LB(I,2)}{\UB(I,1)} \leq \LB(I,2).
\end{align*}
Hence, using the definition of $\LB(I,2)$, we obtain that $f_2(I,x^{\alg})=0$. In particular, $f(I,x^{\alg})$ approximates $(a,0)$ with the desired approximation guarantee in the second component.

\medskip

Now assume that a point with image~$(0,b)$ exists on the Pareto curve for a given instance~$I$. We apply $\alg$ to the corresponding instance of the weighted sum
problem~$\Pi^{\WS}(\gamma)$ for some $0<\gamma<\frac{\LB(I,1)}{\alpha\cdot\UB(I,2)}$. Since $\alg$ is an $\alpha$-approximation algorithm for the weighted sum problem, the resulting
solution~$x^{\alg}=x^{\alg}(\gamma)$ now satisfies
\begin{align*}
  f_1(I,x^{\alg})+ \gamma\cdot f_2(I,x^{\alg}) \leq \alpha\cdot(0 +\gamma\cdot b) = \alpha\cdot \gamma\cdot b. 
\end{align*}
Using the choice of $\gamma$ and that $\gamma\cdot f_2(I,x^{\alg})\geq 0$ and $b\leq \UB(I,2)$, this implies that
\begin{align*}
	f_1(I,x^{\alg})\leq \alpha\cdot \gamma\cdot b < \frac{b\cdot \LB(I,1)}{\UB(I,2)} \leq \LB(I,1).
\end{align*}
Hence, using the definition of $\LB(I,1)$, we obtain that $f_1(I,x^{\alg})=0$. In particular, $f(I,x^{\alg})$ approximates $(0,b)$ with the desired approximation guarantee in the first component.
Moreover, since $f_1(I,x^{\alg})\geq 0$ and $\gamma>0$, we also obtain that $f_2(I,x^{\alg})\leq \alpha\cdot b\leq \alpha\cdot(1+\frac{2}{\epsilon})\cdot b$. Thus, $f(I,x^{\alg})$ also approximates $(0,b)$ with
the desired approximation guarantee in the second component.

\section*{Appendix~$2$}
We identify two problems in the parametric search algorithm presented in Section~$6$ of~\cite{marathe1998networkdesign}.

\medskip

Marathe et al.~\cite{marathe1998networkdesign} consider general bicriteria network design problems. In a generic bicirteria network design problem, one is given an undirected graph~$G=(V,E)$, two positive-valued objective functions~$f_1$ and $f_2$, and a budget~$B>0$ on the value of the first objective function and the task is to find a subgraph from a given subgraph-class that minimizes the second objective function subject to the budget (upper bound) on the first objective function. For example, if the considered subgraph-class consists of all the spanning trees of~$G$ and the two objective values assigned to a spanning tree are the total costs of the edges in the tree with respect to two different sets of edge costs, one obtains the budget-constrained version of the bicriteria minimum spanning tree problem.

\medskip

In Section~$6$ of~\cite{marathe1998networkdesign}, the authors present a bicriteria approximation algorithm that applies whenever the two objectives in the considered bicriteria network design problem are of the same type (e.g., both are total costs of edges computed using two different sets of edge costs as in the bicriteria minimum spanning tree problem). They claim that, when given $\epsilon>0$ and a polynomial-time $\alpha$-approximation algorithm~$\alg$ for the weighted sum problem (which corresponds to the unicriterion problem of minimizing the first objective function~$f_1$ due to the assumption that both objective functions are of the same type), their method yields is a bicriteria $(\alpha\cdot(1+\epsilon),\alpha\cdot(1+\frac{1}{\epsilon}))$-approximation algorithm.

\medskip

\textbf{Parametric Search Algorithm (Marathe et al.~\cite{marathe1998networkdesign})}
\begin{enumerate}
	\item Input: Graph~$G$, budget~$B>0$, accuracy parameter~$\epsilon>0$, $\alpha$-approxi\-mation algorithm~$\alg$ for the weighted sum problem.
	\item Compute an upper bound~$\UB(2)$ on the second objective function value~$f_2(x)$ of all feasible subgraphs~$x$ with $f_1(x)\leq B$.
	\item For any $D>0$, let $h(D)$ denote the objective function value of the subgraph returned by $\alg$ when applied with the weighted sum objective function $\frac{D}{B}\cdot f_1 + f_2$.
	\item Perform a binary search over the interval $[0,\epsilon\cdot \UB(2)]$ to find a $D'\in [0,\epsilon\cdot \UB(2)]$ such that
	\begin{itemize}
		\item $\frac{h(D')}{D'}>\alpha\cdot(1+\epsilon)$,
		\item $\frac{h(D'+1)}{D'+1}\leq\alpha\cdot(1+\epsilon)$.
	\end{itemize}
	\item If the binary search fails to find a valid $D'$, return ``\textit{No Solution}'', else return the solution computed by $\alg$ when applied with the weighted sum objective function $\frac{D'+1}{B}\cdot f_1 + f_2$.
\end{enumerate} 

The correctness and approximation guarantee of the algorithm are shown in~\cite{marathe1998networkdesign} by establishing two claims. The first one (Claim~$6.1$ in~\cite{marathe1998networkdesign}) states that the binary search performed in the algorithm is well-defined. To establish this claim, the authors try to show that the fraction~$\frac{h(D)}{D}$ is monotone nonincreasing in $D$. However, the following example shows that this function is not monotone in general, which implies that the binary search performed in the algorithm is not well-defined:

\begin{example}\label{ex1}
Consider the instance of the budget-constrained version of the bicriteria minimum spanning tree problem given by the budget~$B=2$ and the graph shown in Figure~\ref{fig:ex1}.
\begin{figure}[h]
	\centering
	\begin{tikzpicture}
	\node[circle,minimum width=0.5mm,draw] (a) at (0,0) {1};
	\node[circle,minimum width=0.5mm,draw] (b) at (1.5,2) {2};
	\node[circle,minimum width=0.5mm,draw] (c) at (3,0) {3};
	\draw (a) -- (b) node[pos=0.5, sloped, above]{(3,1)};
	\draw (c) -- (b) node[pos=0.5, sloped, above]{(1,3)};
	\draw (a) -- (c) node[pos=0.5, sloped, below]{(1,1)};
	\end{tikzpicture}
	\caption{Graph with two positive costs on each edge used in Example~\ref{ex1}.}\label{fig:ex1}
\end{figure}
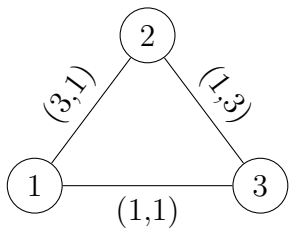

Obviously, there are three possible spanning trees~$x_1$, $x_2$, and $x_3$ with objective function values $(f_1(x_1),f_2(x_1))=(4,2)$, $(f_1(x_2),f_2(x_2))=(2,4)$ and $(f_1(x_3),f_2(x_3))=(4,4)$.
Since we consider a budget of $B=2$, we can use the upper bound $\UB(2)=4$ within the algorithm. The the accuracy parameter is set to~$\epsilon=1$. Although there are exact algorithms that solve the weighted sum problem in polynomial time in this case, suppose that a $\frac{5}{4}$-approximation algorithm~$\alg$ is used. 

We now evaluate the fraction~$\frac{h(D)}{D}$ at two points of the interval $[0,4]=[0,\epsilon\cdot \UB(2)]$. To do so, let $\hat{h}(x,D)\colonequals \frac{D}{B}\cdot f_1(x)+f_2(x)$ denote the objective function value of spanning tree~$x$ in the weighted sum problem for parameter~$D$.
\begin{itemize}
	\item For $D_1\colonequals 3$, we get $\hat{h}(x_1,3)=8$, $\hat{h}(x_2,3)=7$, and $\hat{h}(x_3,3)=10$. Here, $x_2$ is the optimal solution. If $\alg$ returns this solution, we get $\frac{h(3)}{3}=\frac{7}{3}$.
	\item For $D_2\colonequals 4$, we get $\hat{h}(x_1,4)=10$, $\hat{h}(x_2,4)=8$, and $\hat{h}(x_3,4)=12$. Although~$x_2$ is again the optimal solution, it might happen that $\alg$ returns~$x_1$ (which does not contradict the approximation guarantee of $\alg$ since $10\leq \frac{5}{4}\cdot 8$). Then we obtain $\frac{h(4)}{4}=\frac{5}{2}$. 
\end{itemize} 
Consequently if $\alg$ returns the solutions as stated above, we obtain that $\frac{h(3)}{3}=\frac{7}{3}<\frac{5}{2}=\frac{h(4)}{4}$ and $\frac{h(D)}{D}$ is \emph{not} monotone nonincreasing in~$D$ on the considered interval $[0,\epsilon\cdot \UB(2)]=[0,4]$.
\end{example}

The reason that the claimed monotonicity of the fraction~$\frac{h(D)}{D}$ fails is that an approximation algorithm is used for the weighted sum problem. If an \emph{exact} algorithm is used for the weighted sum problem, the monotonicity follows from the proof of Claim~$6.1$ given in~\cite{marathe1998networkdesign}. In the general case where $\alg$ is only an approximation algorithm, however, the algorithm cannot be used due to the ill-defined binary search (this problem cannot be circumvented by replacing the binary search by a linear search since this would not yield a polynomial-time algorithm anymore).

\medskip

We now provide a second example showing that, even if an exact algorithm is used to solve the weighted sum problem, Claim~$6.2$ in~\cite{marathe1998networkdesign} (which states that the algorithm always outputs a solution obtaining the desired approximation guarantee) does not always hold:

\begin{example}\label{ex2}
Consider the instance of the budget-constrained version of the bicriteria minimum spanning tree problem given by the budget~$B=3$ and the graph shown in Figure~\ref{fig:ex2}. 
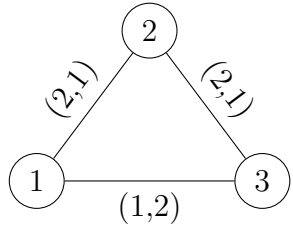
\begin{figure}[h]
	\centering
	\begin{tikzpicture}
		\node[circle,minimum width=0.5mm,draw] (a) at (0,0) {1};
		\node[circle,minimum width=0.5mm,draw] (b) at (1.5,2) {2};
		\node[circle,minimum width=0.5mm,draw] (c) at (3,0) {3};
		\draw (a) -- (b) node[pos=0.5, sloped, above]{(2,1)};
		\draw (c) -- (b) node[pos=0.5, sloped, above]{(2,1)};
		\draw (a) -- (c) node[pos=0.5, sloped, below]{(1,2)};
	\end{tikzpicture}
	\caption{Graph with two positive costs on each edge used in Example~\ref{ex2}.}\label{fig:ex2}
\end{figure}
Again, there are three possible spanning trees~$x_1$, $x_2$, and $x_3$. The corresponding objective function values are $(f_1(x_1),f_2(x_1))=(4,2)$, and $(f_1(x_2),f_2(x_2))=(f_1(x_3),f_2(x_3))=(3,3)$. Since we consider a budget of $B=3$, we can use the upper bound~$\UB(2)=3$ within the algorithm. The accuracy parameter is set to $\epsilon=\frac{2}{3}$, so we obtain $[0,\epsilon\cdot \UB(2)]=[0, 2]$.
In order to solve the weighted sum problem (which is a unicriterion minimum spanning tree problem), we use an exact algorithm~$\alg$ (so $\alpha=1$).

\medskip

\noindent
Now consider the values~$h(D)$ for $0\leq D\leq 3$. For those values of~$D$, we have (using the notation introduced in Example~\ref{ex1})
\begin{align*}
	\hat{h}(x_1,D)=\frac{4}{3}\cdot D+2\leq D+3=\hat{h}(x_2,D)=\hat{h}(x_3,D),
\end{align*}
so $h(D)=\frac{4}{3}\cdot D+2$. Hence, for all values of $D$ in the considered interval $[0,\epsilon\cdot \UB(2)]=[0, 2]$, we obtain
\begin{align*}
	\frac{h(D+1)}{D+1}=\frac{4}{3}+\frac{2}{D+1}\overset{D\leq 2}{\geq}\frac{4}{3}+\frac{2}{3}>1+\frac{2}{3}=\alpha\cdot(1+\epsilon).
\end{align*} 
Consequently, there is no $D\in[0,2]=[0,\epsilon\cdot \UB(2)]$ that fulfills both inequalities in step~$4$ of the algorithm and the algorithm returns ``\textit{No Solution}''. However, with $x_2$ and $x_3$, there exist two feasible solutions that satisfy the given budget of~$B=3$ on the first objective function. 
\end{example}
\end{document}